\documentclass[10pt]{article}

\usepackage{graphics}
\usepackage{color}
\usepackage{bm}
\usepackage{amsfonts}

\input epsf

\begin{document}

\title{
Daehee Formula Associated with the $q$-Extensions of Trigonometric
Functions
\thanks{This work was partially  supported by Jangjeon Mathematical Society(JMS2006-12-C0007).
} }

\author{Taekyun Kim,
        EECS, Kyungpook National University, \\
        Taegu, 702-701. South Korea \\
e-mail: tkim@knu.ac.kr\\\\
        Soonchul Park,  EECS, Kyungpook National University, \\
        Taegu, 702-701. South         Korea\\
        e-mail:scp@knu.ac.kr\\\\
        Seog-Hoon Rim, Department of Mathematics Education,\\
        Kyungpook National University, Taegu, 702-701, S. Korea\\
        e-mail: shrim@knu.ac.kr
        }


\maketitle

\begin{abstract}
In this paper we introduce a Daehee constant which is called
$q$-extension of Napier constant, and consider Daehee formula
associated with the q-extensions of trigonometric functions.
 That is, we derive the $q$-extensions of sine and cosine
functions from this Daehee formula. Finally, we give the
$q$-calculus related to the $q$-extensions of sine and cosine
functions.

\end{abstract}

\smallskip


\smallskip

{\bf AMS(MOS) subject classifications.} 11S80, 11B98


\section{ Introduction}

When one talks of $q$-extension, $q$ is variously considered as an
indeterminate, a complex number $q \in \mathbb{C}$, or a $p$-adic
number $q \in \mathbb{C}_p$. The $q$-basic numbers are defined
\cite{Kim, Kim Rim, Koornwinder, Schork, Simsek} by
\[
[n]_q = \frac{1-q^n}{1-q} = 1 + q + q^2 + \cdots + q^{n-1},
\]
and the $q$-factorials by $[n]_q! = [n]_q [n-1]_q \cdots [1]_q$.
Newton's binomial formula says:
\begin{equation}\label{Newton binomial formula}
(x +y)^n = \sum_{k = 0}^n \left(\begin{array}{c} n \\ k \end{array}
\right) y^{n-k} x^k, \;\; n \in \mathbb{Z}_{+},
\end{equation}
where
\begin{equation}\label{n comb k}
\left(\begin{array}{c} n \\ k \end{array} \right) = \frac{n (n-1)
\cdots (n - k + 1)}{k!} = \frac{n!}{k! (n-k)!},
\end{equation}
which is called binomial coefficient.

Throughout this paper, let us assume that $q \in \mathbb{C}$ with
$|q| < 1$. Here we use the $q$-deformed binomial coefficient (or
Gaussian binomial coefficient) which is also defined by
\begin{equation}\label{q-deformed binomial coefficients}
\left(\begin{array}{c} n \\ k \end{array} \right)_q =
\frac{(q:q)_n}{(q:q)_k (q:q)_{n-k}} = \frac{[n]_q !}{[n-k]_q!
[k]_q!},
\end{equation}
while the $q$-shifted factorial is given by
\[
(a:q)_k = (1-a)(1-aq) \cdots (1-aq^{k-1}), \; \; a \in \mathbb{C},
\; k \in \mathbb{Z}_{+}.
\]
The recurrence relations below show that the $q$-binomial
coefficient is a polynomial in $q$.
$$ \left(\begin{array}{c} x \\ k \end{array} \right)_q = q^k \left(\begin{array}{c} x-1 \\ k \end{array} \right)_q
+\left(\begin{array}{c} x-1 \\ k-1 \end{array} \right)_q  .$$

Note that the $q \rightarrow 1$ yields the conventional number
$[n]_{q=1} = n$ and therefore, the conventional binomial
coefficient (see \cite{Kim Rim, Koornwinder, Schork});
\[
\left(\begin{array}{c} n \\ k \end{array} \right)_{q=1} =
\left(\begin{array}{c} n \\ k \end{array} \right).
\]

It is well known that
\begin{equation}\label{Napier constant}
\lim_{x \rightarrow \infty} \left( 1 + \frac{1}{x} \right)^x =
\lim_{x \rightarrow 0} \left( 1 + x \right)^{\frac{1}{x}} = e,
\end{equation}
which is called Napier constant (or Euler number). Let $i$ be a
complex that is defined by $i = \sqrt{-1} = (-1)^{\frac{1}{2}}$.
Then Euler formula is given by
\begin{equation}\label{Euler formula}
e^{ix} = \cos x + i \sin x, \;\; \mbox{ (see \cite{Kreyszig,
Kim-1}).}
\end{equation}
From ($\ref{Euler formula}$), we derive
\begin{equation}\label{sine and cosine formula1}
\cos x = \frac{1}{2}(e^{ix} + e^{-ix}), \;\; \mbox{ and } \;\;
\sin x = \frac{1}{2i}(e^{ix} - e^{-ix}).
\end{equation}
By using Taylor expansion, we easily see (\cite{Kreyszig}) that
\begin{equation}\label{sine and cosine formula2}
\sin x = \sum_{n=0}^{\infty} \frac{(-1)^n}{(2n+1)!} x^{2n+1}, \;\;
\mbox{ and } \;\; \cos x = \sum_{n=0}^{\infty} \frac{(-1)^n}{(2n)!}
x^{2n}.
\end{equation}

In the recent paper, M. Schork has studied Ward's ``Calculus of
Sequences" and introduced $q$-addition $x \oplus_q y$ in
\cite{Schork} by
\begin{equation}\label{q-addition}
(x \oplus_q y)^n = \sum_{k = 0}^n \left(\begin{array}{c} n \\ k
\end{array} \right)_q x^k y^{n-k}.
\end{equation}
This $q$-addition was already known to Jackson and was generalized
later on by Ward and AI-Salam \cite{Schork}. The two
$q$-exponentials are defined in \cite{Koornwinder} by
\begin{equation}\label{q-exponential1}
e_q(z) = \sum_{n=0}^{\infty}\frac{z^n}{[n]_q!} = \sum_{n=0}^{\infty}
\frac{((1-q)z)^n}{(q:q)_n} = \frac{1}{(z(1-q):q)_{\infty}},
\end{equation}
\begin{equation}\label{q-exponential2}
E_q(z) = \sum_{n=0}^{\infty}\frac{q^{\tiny (\begin{array}{c} n \\
2 \end{array} )}{z^n}}{[n]_q!} = \sum_{n=0}^{\infty}
\frac{((1-q)z)^n q^{\tiny (\begin{array}{c} n
\\ 2 \end{array} )}}{(q:q)_n} = (-z(1-q):q)_{\infty},
\end{equation}
where $(a:q)_{\infty} = \lim_{k \rightarrow \infty} (a:q)_k =
\Pi_{i=1}^{\infty}(1-aq^{i-1})$ (see \cite{Koornwinder, Schork,
Kim-1}). From ($\ref{q-exponential1}$) and
($\ref{q-exponential2}$), we note that $e_q(z) \cdot E_q(-z) = 1$
for $|z| < 1$.

The purpose of this paper is to construct Daehee formula and the
$q$-extensions of sine and cosine. From these $q$-extensions, we
derive some interesting formulae related to Daehee formula and the
$q$-extension of sine and cosine functions.

\section{$q$-extension of Euler formula and trigonometric functions}

Let us consider the Jackson $q$-derivative $D_q$ by
\begin{equation}\label{Jackson q-derivative}
D_q f(x) = \frac{f(x) - f(qx)}{(1-q)x}, \;\; \mbox{ see \cite{Kim
Rim, Koornwinder, Schork}.}
\end{equation}
 It satisfies $D_q x^n =
[n]_q x^{n-1}$ and reduces in the limit $q \rightarrow 1$ to the
ordinary derivative. From the definition of $q$-exponential
function, we derive
\[
D_q(e_q(\lambda x)) = \lambda e_q(\lambda x),
\]
since
\[
D_q\left( \sum_{n=0}^{\infty} \frac{\lambda^n x^n}{[n]_q!} \right) =
\sum_{n=1}^{\infty} \frac{\lambda^n x^{n-1}}{[n-1]_q!} =
\sum_{n=0}^{\infty} \frac{\lambda^{n+1} x^n}{[n]_q!} = \lambda e_q
(\lambda x), \;\; \mbox{ see \cite{Schork}. }
\]
The $q$-integral was defined in \cite{Kim Rim} by
\begin{equation}\label{q-integral}
\int_0^x f(t) d_qt = (1-q) \sum_{k=0}^{\infty} f(q^k x) q^k x,
\end{equation}
where $x \in \mathbb{R}$, and the right-hand side converges
absolutely. In particular, if we take $f(x) = x^n$, then we have
\[
\int_0^x t^n d_q t = (1-q) \sum_{k=0}^{\infty} q^{(n+1)k} x^{n+1} =
\frac{x^{n+1}}{[n+1]_q}, \;\; \mbox{ see \cite{Kim Rim}. }
\]
From ($\ref{Jackson q-derivative}$) and ($\ref{q-integral}$), we
derive
\[
\int_0^x D_q f(t) d_q t = \sum_{k=0}^{\infty}(f(q^k x) -
f(q^{k+1}x)) = f(x), \;\; \mbox{ cf. \cite{Kim Rim}. }
\]
Therefore we have the following:

\smallskip

{\bf Lemma 1}. \quad {\it Let $f$ be a $q$-integrable function. Then
we have
\[
\int_0^x D_q f(t) d_qt = f(x).
\]
}

\smallskip

By the definition of $q$-integral, we easily see that
\[
D_q(f(x)g(x)) = f(x) D_q g(x) + g(qx) D_q f(x),
\]
\[
\int_0^x f(t) D_q(g(t)) d_qt = f(x)g(x) - \int_0^x g(qt)
D_q(f(t))d_qt.
\]

Let us consider the $q$-extension of Napier constant (or Euler
number) which is called Daehee constaint:
\begin{eqnarray}
\lim_{x\rightarrow \infty} \left( 1 \oplus_q \frac{1}{[x]_q}
\right)^x & = & \lim_{x\rightarrow \infty} \sum_{k=0}^{\infty}
\left(\begin{array}{c} x \\ k \end{array} \right)_q \left(
\frac{1}{[x]_q}\right)^k \nonumber \\
& = & \sum_{k=0}^{\infty} \lim_{x \rightarrow \infty} \frac{(1-q^x)
\cdots (1-q^{x-k+1})}{[k]_q! (1-q^x)^k} \nonumber \\
& =& \sum_{k=0}^{\infty} \frac{1}{[k]_q!} = e_q.
\end{eqnarray}

\smallskip

{\bf Proposition 2}. \quad {( Daehee constant)}
\[
e_q = \lim_{x \rightarrow \infty} \left( 1 \oplus_q \frac{1}{[x]_q}
\right)^x.
\]

\smallskip

From the definition of $e_q(x)$, we consider
\begin{equation}\label{eqix}
e_q(ix) = \sum_{n=0}^{\infty} \frac{(ix)^n}{[n]_q!} =
\sum_{n=0}^{\infty}\frac{(-1)^n x^{2n}}{[2n]_q!} + i
\sum_{n=0}^{\infty}\frac{(-1)^n x^{2n+1}}{[2n+1]_q!}
\end{equation}
In the viewpoint of (\ref{sine and cosine formula2}), we define the
$q$-extension of sine and cosine function as follows:
\begin{equation}\label{q-ext sine and  cosine}
\sin_q x = \sum_{n=0}^{\infty}\frac{(-1)^n x^{2n+1}}{[2n+1]_q!},
\;\; \mbox{ and } \;\; \cos_q x = \sum_{n=0}^{\infty}\frac{(-1)^n
x^{2n}}{[2n]_q!}
\end{equation}
From ($\ref{eqix}$) and ($\ref{q-ext sine and  cosine}$), we
derive Daehee formula which is the $q$-extension of Euler formula:
\[
e_q(ix) = \cos_q x + i \sin_q x.
\]
Note that
\[
\lim_{q \rightarrow 1} e_q(ix) = \cos x + i \sin x.
\]
Therefore we obtain the following:

\smallskip
\smallskip

{\bf Theorem 3}. \quad {(Daehee formula)}
\[
e_q(ix) = \cos_q x + i \sin_q x.
\]

\smallskip

By the definition of $q$-addition, we easily see that
\begin{eqnarray*}
e_q(x) e_q(y) &=& \left( \sum_{k=0}^{\infty} \frac{x^k}{[k]_q!}
\right) \left( \sum_{l=0}^{\infty} \frac{y^l}{[l]_q!} \right)
 = \sum_{n=0}^{\infty} \left( \sum_{k=0}^{n} \frac{x^k y^{n-l}}{[k]_q!
[n-k]_q!}\right) \\
 &=& \sum_{n=0}^{\infty} \left( \sum_{k=0}^{n} \frac{[n]_q!}{[k]_q!
[n-k]_q!}x^k y^{n-l}\right) \frac{1}{[n]_q!}
 = \sum_{n=0}^{\infty} \left( \sum_{k=0}^{n} \left(\begin{array}{c} n
\\ k \end{array} \right)_qx^k y^{n-l}\right) \frac{1}{[n]_q!} \\
 &=& \sum_{k=0}^{\infty}\frac{(x \oplus_q y)^n}{[n]_q!} = e_q(x \oplus_q y).
 \end{eqnarray*}
 Therefore we obtain the following:

\smallskip

 {\bf Proposition 4}. \quad {\it For $x, y \in \mathbb{R}$, we have
 \[
 e_q (x) e_q (y) = e_q (x \oplus_q y).
 \]
}

\smallskip

 From proposition 4, we can derive the following:
 \[
 e_q(ix) e_q(iy) = e_q(i(x \oplus_q y)).
 \]
 By using Daehee formula, we see that
 \[
 (\cos_q x + i\sin_qx)(\cos_qy + i\sin_qy) = \cos_q(x\oplus_qy) + i
 \sin_q(x \oplus_qy).
 \]
 Thus, we have
 \[
 \begin{array}{l}
 (\cos_qx\cos_qy - \sin_qx\sin_qy) + i(\sin_qx\cos_qy + \cos_qx\sin_qy) \\
 = \cos_q(x\oplus_qy) + i\sin_q(x\oplus_qy).
 \end{array}
 \]
 By comparing the coefficients on both sides, we obtain the
 following:

\smallskip

 {\bf Theorem 5}. \quad {\it For $x, y \in \mathbb{R}, \; q \in
 \mathbb{C}$ with $|q| < 1$, we have
\begin{eqnarray*}
\cos_q(x\oplus_qy) = \cos_qx\cos_qy - \sin_qx\sin_qy, \\
\sin_q(x\oplus_qy) = \sin_qx\cos_qy + \cos_qx\sin_qy.
\end{eqnarray*}
}

\smallskip

By the same motivation of ($\ref{q-addition}$), we can also define
\begin{equation}\label{q-subtration}
(x\ominus_qy)^n = \sum_{k=0}^n \left(\begin{array}{c} n
\\ k \end{array} \right)_q (-1)^{n-k}x^ky^{n-k}.
\end{equation}
By ($\ref{q-ext sine and  cosine}$), we easily see that
\begin{equation}\label{odd and even property sine and cosine}
\sin_q(-x) = -\sin_q(x) \;\; \mbox{ and } \cos_q(-x) = \cos_q(x).
\end{equation}
From Theorem 5, ($\ref{q-subtration}$) and ($\ref{odd and even
property sine and cosine}$), we note that
\begin{eqnarray*}
\cos_q^2x + \sin_q^2x &=& \cos_q(x(1\ominus_q 1)), \\
\cos_q(x(1\oplus_q1)) &=& \cos_q^2x - \sin_q^2 x, \\
\sin_q(x(1\oplus_q1)) &=& 2\sin_qx \cos_qx.
\end{eqnarray*}
Therefore we obtain the following:

\smallskip

{\bf Corollary 6}. {\it For $x, y \in \mathbb{R}$, we have
\begin{eqnarray*}
\cos_q^2x + \sin_q^2x &=& \cos_q(x(1\ominus_q 1)), \\
\cos_q^2x - \sin_q^2x &=& \cos_q((1\oplus_q1)x), \\
2\sin_qx \cos_qx &=& \sin_q(x(1\oplus_q1)).
\end{eqnarray*}
}

\smallskip

By using Daehee formula, we easily see that
\begin{eqnarray*}
e_q(ix) = \cos_qx + i\sin_qx, \\
e_q(-ix) = \cos_qx - i\sin_qx.
\end{eqnarray*}
Thus, we have the following:

\smallskip

{\bf Corollary 7}. {\it For $x, y \in \mathbb{R}$, we have
\[
\cos_qx = \frac{e_q(ix) + e_q(-ix)}{2}, \;\; \mbox{ and } \; \;
\sin_qx = \frac{e_q(ix) - e_q(-ix)}{2i}
\]
}

\smallskip

We now define $q$-extension of $\tan x$, $\sec x$, $\csc x$ and
$\cot x$ as follows:
\[
\tan_qx = \frac{\sin_qx}{\cos_qx}, \; \sec_qx = \frac{1}{\cos_qx},
\; \csc_qx = \frac{1}{\sin_qx}, \mbox{ and } \cot_qx =
\frac{\cos_qx}{\sin_qx}.
\]
Note that
\begin{equation}\label{properties of tan and cot}
\begin{array}{lllll}
1 + \tan_q^2x &=& \frac{\sin_q^2x + \cos_q^2x}{\cos_q^2x} &=&
\cos_q(x(1\ominus_q 1))\sec_q^2x, \\
1 + \cot_q^2x &=& \frac{\cos_q^2x + \sin_q^2x}{\sin_q^2x} &=&
\cos_q(x(1\ominus_q 1))\csc_q^2x.
\end{array}
\end{equation}
From ($\ref{Jackson q-derivative}$), we can derive
\begin{equation}\label{q-derivative of rat fn}
D_q\left(\frac{f(x)}{g(x)}\right) = \frac{1}{g(x)g(qx)}
[D_q(f)g(qx) - f(qx)D_q(g)].
\end{equation}
By ($\ref{Jackson q-derivative}$), ($\ref{q-ext sine and  cosine}$)
and ($\ref{q-derivative of rat fn}$), we easily see that
\[
D_q(\sin_qx) = \cos_qx, \;\; \mbox{ and } \;\; D_q(\cos_qx) =
-\sin_qx.
\]
Moreover,
\[
D_q(\tan_qx) = 1 + \tan_qx \tan_q(qx).
\]
Therefore we obtain the following.

\smallskip

{\bf Theorem 8}. \quad {\it For $x \in \mathbb{R}$, we have
\[
D_q(\sin_qx) = \cos_qx, \;\; \mbox{ and } \;\; D_q(\cos_qx) =
-\sin_qx.
\]
Moreover,}
\[
D_q(\tan_qx) = 1 + \tan_qx \tan_q(qx).
\]

\smallskip

By Lemma 1 and Theorem 8, we have the following:

\smallskip

{\bf Corollary 9}. \quad {\it For $x \in \mathbb{R}$, we have
\[
\int_0^x \sin_qt d_qt = -\cos_qx, \;\; \mbox{ and } \;\; \int_0^x
\cos_qt d_qt = -\sin_qx.
\]
}

\smallskip

Finally,
\[
\int_0^x (1 + \tan_qt \tan_q(qt)) d_qt = \tan_qx.
\]

\end{document}